\documentclass[12pt]{article}
\def\part#1{\frac{\partial\phantom{q}}{\partial#1}}
\newtheorem{thm}{Theorem}
\newtheorem{prp}[thm]{Proposition}

{\hfill $\Box$ \end{trivlist}}
{\end{trivlist}}
\newenvironment{rmk}{\begin{trivlist}\item[]{\bf Remarks} }%
{\end{trivlist}}

\def\vol{\mathop{\rm vol}\nolimits}

\def\sgn{\mathop{\rm sgn}\nolimits}

\def\tr{\mathop{\rm tr}\nolimits}
\def\ch{\mathop{\rm ch}\nolimits}
\def\Td{\mathop{\rm Td}\nolimits}
\def\End{\mathop{\rm End}\nolimits}

\newcommand{\R}{\mathbf{R}}

\newcommand{\Z}{\mathbf{Z}}

\begin{document}
\title{Curvature and characteristic numbers of  hyperk\"ahler manifolds}
\author{Nigel Hitchin and  Justin Sawon\footnote{On leave from Trinity College, Cambridge}\\[5pt]
\itshape Mathematical Institute\\
\itshape 24-29 St Giles\\
\itshape Oxford OX1 3LB\\
\itshape England\\[12pt]}
\maketitle
\section{Introduction}
In this paper we prove the following formula for the ${\cal L}^2$ norm of the Riemann curvature tensor $R$ of an irreducible compact hyperk\"ahler manifold $M$ of real dimension $4k$:
\begin{equation}
\frac{1}{(192\pi^2 k)^k}  \frac{\Vert R \Vert ^{2k}}{(\vol M)^{k-1}}=\sqrt{\hat A}\,[M]
\label{A}
\end{equation}
By $\sqrt{\hat A}$ we mean here the multiplicative sequence of Pontryagin classes defined by the power series
$$\left(\frac{{\sqrt z}/2}{\sinh {\sqrt z}/2}\right)^{1/2}$$
It is perhaps surprising to see this norm expressed solely in terms of the volume and Pontryagin classes. In fact the context for this formula is the work of Rozansky and Witten, who  showed in \cite{RW} how one could associate to each
compact hyperk\"ahler manifold  a three-manifold invariant. In more abstract 
terms,
what they show is that $M^{4k}$ defines a class in the dual of the graph
homology group ${\cal A}^k(\phi)$ consisting of linear combinations of oriented
trivalent graphs with $2k$ vertices, modulo antisymmetry and the IHX relation. To evaluate the class on a
given graph, one uses the graph to 
 perform contractions on the $2k$-th tensor power of the curvature tensor and then
 integrates over $M^{4k}$ the resulting scalar.
 
 Characteristic numbers (in the hyperk\"ahler case these are just Pontryagin numbers)  are obtained by
 integrating polynomials in the curvature, and these also fall within the Rozansky-Witten approach. What we do in this paper is to  organize the evaluation of characteristic numbers using  the  graph-homological formalism. The benefits of this are firstly that the IHX relation offers a systematic way of performing integration by parts to relate invariants corresponding to different graphs, and secondly that the graphs involved have fortuitously been studied  in 3-manifold theory. We make use of the recently proved Wheeling Theorem  \cite{blt} to obtain our result. This enables us to show that the graph consisting of $k$ disjoint copies of the unique connected trivalent graph with two vertices (which we denote by $\Theta$ for obvious reasons)  is homologous to  a specific linear combination of 
 graphs defining characteristic numbers. It is not hard to see that the Rozansky-Witten invariant for $\Theta^k$ is expressible directly in terms of the norm of the curvature and the volume, and putting the two together we obtain the formula (\ref{A}). 

An immediate corollary of the formula is the inequality (for compact irreducible hyperk\"ahler manifolds)
$$\sqrt{\hat A}\,[M]>0.$$
The  authors wish to thank J. Ellegaard Andersen and S. Willerton for useful
 conversations. They are also   grateful  to  S. Garoufalidis  for  enlightening  them on  many aspects of  the  subject  of graph homology.  

 \section{Rozansky-Witten invariants}
 We  review here for the reader's benefit the construction of Rozansky and
 Witten \cite{RW},  and  also  its   interpretation by Kapranov using Dolbeault
 cohomology \cite{Kap}.   
  Recall that a hyperk\"ahler manifold is a Riemannian manifold $M^{4k}$ with
  holonomy   contained in $Sp(k)$. We shall say that $M$ is irreducible if the holonomy is precisely $Sp(k)$. The metric is k\"ahlerian with respect to
  integrable complex structures $I,J,K$ which act as the quaternions on the
  tangent bundle.   The   corresponding   K\"ahler forms are $\omega_1,\omega_2,\omega_3$. 
  If we fix attention on the
  complex structure $I$, then 
  $\omega=\omega_2+i\omega_3$
  is a covariant constant holomorphic symplectic form. 
  
  The Riemann curvature
  tensor $R$ can be considered, relative to the complex structure $I$, as a section   $K\in \Omega^{1,1}(\End T)$ with components $K^i_{jk\bar l}$ relative  to local complex coordinates.   Using  the non-degenerate holomorphic
  2-form      $\omega$   to   identify   $T$   and   $T^*$, we can lower the
  first 
  index and define  $$\Phi \in \Omega^{1,1}(T^*\otimes T^*)=\Omega^{0,1}(T^*\otimes T^*\otimes T^*)$$ by  $$\Phi_{ijk\bar l}=\sum_m\omega_{im}K^m_{jk\bar l}$$   This is symmetrical in $j,k$ because the connection is torsion-free and
  preserves the complex structure. It is also symmetrical in $i,j$ because the
  the curvature takes values in the Lie algebra of
  $Sp(2k, {\bf C})$, and this consists of matrices of the form $A^i_j$ where 
  $S_{ij}=\sum_k\omega_{ik}A^k_j$ with $S_{ij}$ symmetric.
     Thus   $$\Phi \in \Omega^{0,1}(Sym^3T^*)$$
  The Bianchi identity tells us that $\bar \partial \Phi=0$, and so $\Phi$
  defines a  Dolbeault cohomology class 
  $$[\Phi]\in H^1(M, Sym^3T^*)$$
  \vskip .25cm
  Now let $\Gamma$ be a trivalent graph with $2k$ vertices and no edges joining a vertex to itself. One defines an
  orientation on $\Gamma$ to be an equivalence class of cyclic orderings at
  each vertex. Two such orderings are equivalent if they 
  differ   on an even number of vertices. This definition is good for visualizing
  graphs -- drawing a graph in the plane gives it an induced orientation. An
  equivalent formulation which is more directly applicable to the
  Rozansky-Witten 
  invariants   defines an orientation to be an equivalence class of orderings 
  of vertices   together with orientations on the  edges. If the orderings
  differ   by a permutation $\pi$ and the orientations differ on $n$ edges, then
  the pairs are equivalent if $\sgn \pi=(-1)^n$. We shall see that these two definitions are equivalent in Section 4.
  \vskip .25cm
  Choose an ordering of the vertices and consider the tensor $\Phi\otimes
  \Phi \otimes\dots\otimes\Phi$ with $2k$ factors. If vertex $v_m$ and vertex
  $v_n$   of   the graph $\Gamma$ 
  are   joined by an edge and $m< n$ then we contract with the skew form $\tilde \omega$ on $T^*$ dual to 
  $\omega$ (note in passing that in a dual basis the matrix  $\omega^{ij}$ of $\tilde \omega$ is  the {\it negative} of the inverse of $\omega_{ij}$)
  $$c_{m,n}\omega^{i_m i_n}\Phi\otimes \dots \otimes\Phi_{i_m\dots}\otimes
  \dots   \otimes\Phi_{i_n\dots}\otimes \dots \otimes \Phi$$
  where $c_{m,n}=1$ if the orientation on the edge goes from $v_m$ to $v_n$
  and $c_{m,n}=-1$ if it is in the opposite direction. Continuing over all
  $3k$ edges, we are left with a section of $\bar T^*\otimes \dots \otimes
  \bar T^*$. Projecting this to the exterior product, we obtain
  $$\Gamma (\Phi)\in \Omega^{0,2k}$$
  It is clear now how the orientation on the graph enters, since a change of
  orientation on an edge produces a change in sign of $\Gamma(\Phi)$ since
  $\tilde \omega$ is skew-symmetric, and a change in the sign of the permutation of
  the vertices changes the sign of the exterior product.
  \vskip .25cm
  The Rozansky-Witten invariant of $M^{4k}$ defined by the graph
  $\Gamma$   is      
  \begin{equation}
  b_{\Gamma}(M)=\frac{1}{(8\pi^2)^{k}k!}\int_M \Gamma(\Phi)\omega^k
  \label{def}
  \end{equation}
  
  In Dolbeault terms, we take the class $[\Phi]\in H^1(M, Sym^3T^*)$ and
  perform the same contractions with the holomorphic form $\tilde \omega$ to
  obtain a class in $H^{2k}(M, {\cal O})$ represented by the $(0,2k)$-form
  $\Gamma(\Phi)$. Since $\omega^k\in H^0(M,\Omega^{2k})$, $b_{\Gamma}(M)$ is essentially 
  the   Serre duality pairing of $\omega^k$ and $[\Gamma(\Phi)]$.
  \vskip .15cm
  \begin{rmk}
 
  \item 
  
  \noindent 1. The definition of the invariant given here seems to involve a
  choice   $I$ of complex structure. In fact this is not the case. The best way to see
  this is to think of a hyperk\"ahler structure as in \cite{Sal} as an
  isomorphism $T^c\cong E\otimes {\bf C}^2$, where $E$ is a certain complex
  vector bundle with a non-degenerate skew form. In four dimensions it is the 
  anti-self-dual   spinor   bundle.   The curvature is then a section of $Sym^4 E$ and the 
  contractions on this give the integrand (as in Rozansky and Witten's original paper \cite{RW}).
  This is then clearly independent of  $I,J$ or $K$.
  A particular  choice   of   complex   structure   identifies $E$ with the
  complex   tangent bundle, and this is what we have done above. Using this
  formalism and the quaternionic structure on $E$ it is also easy to see that
  the   invariant $b_{\Gamma}(M)$ is real.
  \item
   \noindent 2. Both the contraction and integration  make use of  the form
  $\omega$. If we multiply $\omega$ by $\lambda$, the invariant $b_{\Gamma}(M)$
 is in fact    unchanged. To see this, note that the curvature is scale
  invariant,   and the definition of $\Phi$     required   a   contraction of the curvature
   with $\omega$, so $\Phi\otimes \dots \otimes
  \Phi$ is scaled by $\lambda^{2k}$. There are $3k$ edges to the graph and to obtain $\Gamma(\Phi)$ we
  contracted with $\tilde \omega$, which scales by $\lambda^{-1}$, once for each edge, so $\Gamma(\Phi)$ is scaled by
  $\lambda^{2k}\lambda^{-3k}=\lambda^{-k}$. This is multiplied by $\omega^k$ to give the final integral
  and since this scales by $\lambda^k$, the result is invariant.
  \item
   \noindent 3. A hyperk\"ahler metric is determined by the three closed forms
  $\omega_1,\omega_2,\omega_3$ \cite{HKLR}. Thus  a first order deformation
  of   a   hyperk\"ahler structure  consists of a triple of closed
  forms   $\dot   \omega_1,   \dot   \omega_2,\dot \omega_3$. These forms are
  of a particular algebraic type -- if the volume is fixed under the deformation then they
  are all of Hodge  type   $(1,1)$ with respect to all   complex structures. Note now
  that 
  a   first   order deformation is  a linear combination of   deformations
     for which we fix two out of the three forms. If we fix   $\omega_2$ and $\omega_3$ then we have fixed the complex  symplectic form
  $\omega$. This determines the complex structure (the $(1,0)$ forms $\alpha$  are   defined by the condition $\alpha \wedge\omega=0$) and the
  holomorphic symplectic form. The remaining variation is   the K\"ahler form. But the  Dolbeault version 
  of the Rozansky-Witten invariant
  shows that it only depends on $M$ as a complex symplectic manifold, thus
  $b_{\Gamma}(M)$ is constant under this deformation. Using the other complex
  structures $J$ and $K$ in turn, we see that $b_{\Gamma}(M)$ is invariant
  under   {\it   any}  first order    deformation.   The Rozansky-Witten invariants are therefore constant on connected components
  of the moduli space of hyperk\"ahler metrics on $M$.
  \item
  \noindent 4. The particular factor in the integral (\ref{def}) is chosen to make the invariant have multiplicative properties.  The product $M^{4k}\times N^{4l}$ of two hyperk\"ahler manifolds is
     again      hyperk\"ahler. If $\Gamma$ is a graph with $2k+2l$ vertices,
     we can form      the invariant $b_{\Gamma}(M\times N)$. To calculate this,   we
     have to make contractions on a tensor product of $2k+2l$ copies of $\Phi
     +\Psi$, with $\Phi \in \Omega^{0,1}(M, Sym^3T^*_M)$ and $\Psi \in
     \Omega^{0,1}(M, Sym^3T^*_N)$ and then take an exterior product on the
     $(0,1)$ factors. Because of this exterior product, any contraction of more
     than $2k$ vertices is zero on $\Phi$. Moreover, any contraction of $\Phi$
     with $\Psi$ is zero. Thus a non-zero contribution involving $\Phi$ can
     only come from a set of $2k$ vertices such that all edges link them. In other words $b_{\Gamma}(M\times N)$ is non-zero only if  $\Gamma$ is the disjoint union of $\gamma$ and 
     $\gamma'$      where      $\gamma$      is      a trivalent graph with
     $2k$ vertices and $\gamma'$ a trivalent graph with      $2l$ vertices. The invariant $b_{\Gamma}(M\times N)$ is then a sum over
     {\it all} such decompositions. If we define the coproduct 
      $$\Delta\Gamma=\sum_{\gamma \sqcup \gamma'=\Gamma} \gamma \otimes
     \gamma'$$
     then with the chosen normalization of (\ref{def})
     \begin{equation}
     b_{\Gamma}(M\times N)= \Delta\Gamma (M,N)=\sum_{\gamma \sqcup
     \gamma'=\Gamma}b_\gamma(M)b_{\gamma'}(N)
     \label{coprod}
     \end{equation}
     with the obvious convention that if the size of $\gamma$ does not match
     the dimension of $M$, the result is zero. Note in particular that if
     $\Gamma$ is 
           connected  then $b_{\Gamma}(M\times N)=0$. 
  \end{rmk}
    \vskip .5cm
    The invariant $b_{\Gamma}(M)$ depends  on the graph through its graph 
    homology     class. To see this, consider the covariant exterior
    derivative  $d_A$ of the Levi-Civita connection $A$ applied to  $\Phi \in \Omega     ^{0,1}(Sym^3T^*)$.  We
    have
    $$d_A^2 \Phi=K(\Phi)\in \Omega^{1,2}(Sym^3T^*)$$
    where $K$ is the curvature. Since $\bar \partial_A \Phi=\bar \partial
    \Phi=0$, this expression can be written as     $$\bar\partial \partial_A
    \Phi=K(\Phi)$$     where $\partial_A \Phi \in \Omega^{1,1}(Sym^3T^*)=\Omega^{0,1}(T^*\otimes
    Sym^3T^*)$. Symmetrizing this gives $S[\partial_A \Phi]\in \Omega^{0,1}(
    Sym^4T^*)$ and 
   $$ \bar\partial (S[\partial_A \Phi])=S[K(\Phi)].$$
    Now the curvature term $S[K(\Phi)]$ is obtained
    by symmetrizing a multiple of the $(0,2)$-form with values in $Sym^4T^*$
    which has components 
    $$ \omega^{ia}\Phi_{ijk\bar l}\Phi_{abc\bar d}$$
    Interchanging $j,k$ or $b,c$ leaves this unchanged by the symmetry of
    $\Phi$. Interchanging the pair $(j,k)$ with $(b,c)$ also leaves it
    invariant because of the skew-symmetry of $\tilde\omega$ and the exterior
    product on the $(0,1)$ terms. Symmetrizing thus involves the sum of three
    terms corresponding to the three different ways of splitting the indices
    $\{j,k,b,c\}$ into two sets of two elements. Thus we can express $
    S[K(\Phi)]$ as a sum of three contractions:
         \begin{equation}
    \bar\partial (S[\partial_A \Phi])=S[K(\Phi)]=C_1(\Phi,\Phi)+C_2
    (\Phi,\Phi)+C_3(\Phi,\Phi)
    \label{IH}
    \end{equation}     
    Now suppose $\Gamma_*$ is an oriented graph with $2k-1$ vertices, whose
    first     vertex     is     4-valent     and     the     rest  trivalent. 
     Contracting     with $\Gamma_*$ gives
    $\Gamma_*(S[\partial_A \Phi],\Phi,\dots ,\Phi)\in \Omega^{2k-1}$  and using
    $\bar     \partial     \Phi=0$ we have from (\ref{IH}),
    $$\bar \partial \Gamma_*(S[\partial_A \Phi],\Phi,\dots
    \Phi)=\Gamma_1(\Phi)+\Gamma_2(\Phi)+\Gamma_3(\Phi)$$
    where $\Gamma_1, \Gamma_2, \Gamma_3$ are the three trivalent graphs
    obtained from the three different ways of expanding  the 4-valent vertex of
    $\Gamma_*$.     Consequently
    $$[\Gamma_1(\Phi)]+[\Gamma_2(\Phi)]+[\Gamma_3(\Phi)]=0\in H^{2k}(M,{\cal
    O})$$
    and 
    $$b_{\Gamma_1}(M)+b_{\Gamma_2}(M)+b_{\Gamma_3}(M)=0$$
    The invariant $b_{\Gamma}(M)$ is thus well defined on the space of linear
    combinations of graphs modulo the IHX relation
    $\Gamma_1+\Gamma_2+\Gamma_3=0$ which is graphically:

\begin{picture}(150,80) (-20,50)

\put(120,100){\line(2,1){10}}
\put(120,100){\line(-2,1){10}}
\put(120,100){\line(0,-1){20}}
\put(120,80){\line(2,-1){10}}
\put(120,80){\line(-2,-1){10}}
\put(135,90){$=$}
\put(155,90){\line(1,0){10}}
\put(155,90){\line(-1,2){8}}
\put(155,90){\line(-1,-2){8}}
\put(165,90){\line(1,2){8}}
\put(165,90){\line(1,-2){8}}
\put(182,90){$-$}
\put(194,106){\line(2,-3){22}}
\put(215,106){\line(-2,-3){22}}
\put(199,82){\line(1,0){11}}

\end{picture}

    Since the invariant is antisymmetric with
    respect to change of orientation, the hyperk\"ahler manifold $M^{4k}$
    defines     a     class in the dual space of the {\it graph homology space}
     ${\cal A}^k(\phi)$ which consists of {\bf Q}-linear combinations of trivalent
     graphs with $2k$ vertices modulo antisymmetry and the IHX relation. This has a product structure given by disjoint union, a coproduct structure from (\ref{coprod}), and is graded by half the number of vertices.
          
     \section{The graph $\Theta^k$}
  
     The invariants from disconnected graphs will be important for us. Their
     special role is simplest to see for irreducible hyperk\"ahler manifolds.
     The fundamental property we shall need for these manifolds is the
     following (see \cite{Bea} for example):
     \vskip .25cm
     \begin{prp}
      Let $M$ be an irreducible compact hyperk\"ahler manifold.
     Then $h^{0,p}=0$ if $p$ is odd and $h^{0,p}=1$ if $p$ is even.
     \end{prp}
     \vskip .25cm
   The result is a standard application of Bochner's vanishing theorem: since     $M$ has
     zero Ricci tensor, all holomorphic $p$-forms are covariant constant and by
     irreducibility, each one is a constant multiple of $\omega^m$. Since
     $h^{p,0}=h^{0,p}$ this means 
     for      the      Dolbeault cohomology groups  that any element
     of $H^{2m}(M,{\cal O})$ is a constant multiple of the class of $\bar \omega^m \in
     \Omega^{0,2m}$.
     \vskip .25cm
     \noindent Now disconnected graphs are disjoint unions of trivalent graphs with fewer than $2k$ vertices. If $\gamma$ has $2m<2k$ vertices, then we can contract using $\tilde \omega$ the class
     $[\Phi]\in H^1(M^{4k},Sym^3T^*)$ to obtain
     $$[\gamma (\Phi)]\in H^{2m}(M, {\cal O})$$
     As we have seen, if $M$ is irreducible, the Dolbeault cohomology class $[\gamma (\Phi)]= c_{\gamma} [\bar
     \omega^m]$ where $c_{\gamma}$ is a constant.  This means in particular
     that if $\Gamma$ has $2k$ vertices and $\Gamma=\gamma \gamma'$ (the product is disjoint union), then
     \begin{equation}
     b_{\Gamma}(M)=\frac{1}{(8\pi^2)^{k}k!}\int_M c_{\gamma}c_{\gamma'}\bar \omega^k\omega^k
\label {b}
\end{equation}
 
We now use the formula for the volume form $V$ of a K\"ahler manifold of complex dimension $n$  with K\"ahler form $\Omega$: 
$V = \Omega^{n}/{n!}$. So, for our hyperk\"ahler manifold, taking the complex structure $J$, $\Omega=\omega_2=(\omega +\bar \omega)/2$ and we have
\begin{equation}
V=\frac{1}{(2k)!}\omega_2^{2k}=\frac{1}{2^{2k}(2k)!}(\omega+\bar\omega)^{2k}=
\frac{1}{2^{2k}(k!)^2}\omega^k\bar\omega^k
\label{vol}
\end{equation}
Hence we can rewrite (\ref{b}) as 
$$b_{\Gamma}(M)=
\frac{k!}{(2\pi^{2})^k}c_{\gamma}c_{\gamma'}\vol(M)$$
       In particular, if $\Gamma=\Theta^k$ then
       \begin{equation}
b_{\Theta^k}(M)=\frac{k!}{(2\pi^{2})^k}c_{\Theta}^k\vol(M)
\label{bk}
\end{equation}  
 To evaluate the constant $c_{\Theta}$, we shall use an approach which brings in characteristic classes from the start. As in most of this paper, it will be convenient to use a fixed complex structure $I$  to perform calculations, and then switch to a choice-free description when appropriate. Recall that from Chern-Weil theory, if $K\in \Omega^{1,1}(\End T)$ is the curvature of a K\"ahler manifold, then the characteristic class $c_1^2-2c_2$ of the tangent bundle is represented by the closed 4-form
\begin{eqnarray*}
\frac{1}{(2\pi i)^2}\tr K^2&=&
-\frac{1}{4\pi^2}\sum K^i_{jk\bar l}K^j_{im\bar n}dz_k\wedge d\bar z_l\wedge dz_m\wedge d\bar z_n\\&=&\frac{1}{4\pi^2}\sum K^i_{jk\bar l}K^j_{im\bar n}dz_k\wedge dz_m\wedge d\bar z_l\wedge d\bar z_n
\end{eqnarray*}
For a hyperk\"ahler manifold the structure group of $T$ is reduced to $Sp(k)\subset SU(2k)$. In this representation the eigenvalues occur in opposite pairs so all the odd Chern classes vanish, thus the form represents $-2c_2$.   
 Note now that, by linear algebra, if $\alpha=\sum \alpha_{ij} dz_i\wedge dz_j$ is a $(2,0)$-form then
 $$ \alpha \wedge \omega^{k-1}=\frac{1}{2k}(\sum \omega^{ij}\alpha_{ij}) \omega^k$$  
 Hence taking de Rham cohomology classes
 $$c_2[\omega]^{k-1}=-\frac{1}{16\pi^2 k}[\sum \omega^{km}K^i_{jk\bar l}K^j_{im\bar n} d\bar z_l\wedge d\bar z_n][\omega]^k$$
But
$$K^i_{jk\bar l}=-\omega^{ia}\Phi_{ajk\bar l}$$
and so 
$$c_2[\omega]^{k-1}=\frac{1}{16\pi^2 k}\left [\sum \omega^{ai}\omega^{jb}\omega^{kc}\Phi_{ajk\bar l}\Phi_{bic\bar m} d\bar z_l\wedge d\bar z_m\right][\omega]^k=\frac{1}{16\pi^2 k}[\Theta(\Phi)][\omega]^k$$
since the contractions correspond to three edges, each pointing from one vertex to another -- the graph $\Theta$.
Hence
\begin{equation}
c_{\Theta}\int_M \omega^k\bar\omega^k=16\pi^2 k\int_M c_2\omega^{k-1}\bar \omega^{k-1}
\label{cg}
\end{equation}
On the other hand, it is well-known (see \cite{Besse} page 80 for example) that on a Ricci-flat K\"ahler manifold of complex dimension $n$ the ${\cal L}^2$ norm of the curvature can be expressed in terms of $c_2$ and the K\"ahler class:
$$\Vert R\Vert^2=\frac{8\pi^2}{(n-2)!}\int_M c_2\Omega^{n-2}$$
and so in our case,
$$\Vert R\Vert^2=\frac{8\pi^2}{(2k-2)!}\int_M c_2\omega_2^{2k-2}=\frac{8\pi^2}{2^{2k-2}(2k-2)!}\int_M c_2(\omega+ \bar\omega)^{2k-2}$$
but since $c_2$ is of type $(2,2)$ whatever the complex structure,
$$\Vert R\Vert^2=\frac{8\pi^2}{2^{2k-2}(2k-2)!}{2k-2 \choose k-1}\int_M c_2\omega^{k-1} \bar\omega^{k-1}=\frac{8\pi^2}{2^{2k-2}((k-1)!)^2}\int_M c_2\omega^{k-1} \bar\omega^{k-1}$$
and so from (\ref{cg}) and (\ref{vol})
\begin{equation}
c_{\Theta}=\frac{1}{2 k}\frac{\Vert R \Vert^2}{\vol(M)}
\label{ct}
\end{equation}
which finally gives from (\ref{bk})

\begin{equation}
b_{\Theta^k}(M)=\frac{k!}{(4 \pi^2 k)^{k}}\frac{\Vert R \Vert^{2k}}{(\vol M)^{k-1}}
\label{bb}
\end{equation}
For $k=1$, in which case the only irreducible compact hyperk\"ahler manifold is the K3 surface, we have
$$b_{\Theta}(M)=2c_2(M)=48$$
\section{Characteristic numbers}

 In discussing the invariant from the graph $\Theta^k$ we already encountered the Chern-Weil form
$$\frac{1}{(2\pi i)^2}\tr K^2$$
and it is useful to think of all the characteristic classes of the complex tangent bundle as being generated by the  classes
$$s_{2m}=\left[ \frac{1}{(2\pi i)^{2m}}\tr K^{2m}\right]\in H^{4m}(M,\Z)$$
The Chern character in particular is given  by
$$\ch(T)=\sum_m \frac{s_{2m}}{(2m)!}$$
Consider the characteristic number $s_{2k}[M^{4k}]$. This is given by integrating ${(2\pi i)^{-2k}}$ times the form
\begin{equation}
\tr K^{2k}=\sum K^i_{ja\bar b}K^j_{kc\bar d}\dots K^l_{ie\bar f} dz_a\wedge d\bar z_b\wedge dz_c\wedge d\bar z_d\dots dz_e\wedge d\bar z_f
\label{form}
\end{equation}
Since $\Phi_{ijk\bar l}=K^i_{jk\bar l}=-\sum \omega^{ia}\Phi_{ajk\bar l}$, we can  translate the cyclic summation over the first two indices defining the trace of a product of matrices into a contraction using the $2k$ edges around the circle of the graph:

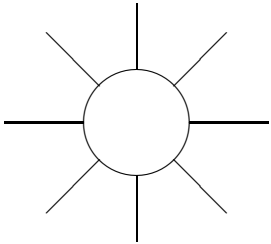
\begin{picture}(150,100) (-100,50)
\put(100,100){\circle{40}}
\put(120,100){\line(1,0){30}}
\put(114,114){\line(1,1){20}}
\put(100,120){\line(0,1){25}}
\put(86,114){\line(-1,1){20}}
\put(80,100){\line(-1,0){30}}
\put(86,86){\line(-1,-1){20}}
\put(100,80){\line(0,-1) {25}}
\put(114,86){\line(1,-1) {20}}
\end{picture}

Such a graph with $2k$ spokes is called a $2k$-{\it wheel} $w_{2k}$. It is not a trivalent graph. The spokes are attached to the hub at trivalent vertices but the end-points of the spokes are univalent. Each spoke interpreted in our evaluation of $\tr K^{2k}$ corresponds to a free $T^*$ index,  the $a$ in $K^i_{ja\bar b}$.

 After contracting around the wheel, $\tr K^{2k}$ now involves taking exterior products of terms like 
$$dz_a\wedge d\bar z_b\wedge dz_c\wedge d\bar z_d\dots \wedge dz_e\wedge d\bar z_f$$
But on reordering this can be written as
\begin{equation}
(-1)^{k}dz_a\wedge dz_c\wedge \dots \wedge dz_e\wedge d\bar z_b\wedge d\bar z_d\dots\wedge d\bar z_f
\label{sign1}
\end{equation}
and this means taking first the exterior product of the free $dz_a$ indices corresponding to the spokes, and following it (as with the general Rozansky-Witten invariant) with the exterior product of the $d\bar z_b$ terms.

Now the algebra of exterior products tells us that
$$\omega^k(v_1,\dots,v_{2k})=\frac{1}{(2k)!}\sum_{\pi \in S_{2k}}\sgn (\pi)\omega(v_{\pi (1)},v_{\pi (2)})\omega(v_{\pi (3)},v_{\pi (4)})\dots \omega(v_{\pi (2k-1)},v_{\pi (2k)})$$
and dually we find
\begin{equation}
\theta_1\wedge \dots\wedge \theta_{2k}=\frac{1}{2^{2k}(k!)^2}\sum_{\pi \in S_{2k}}\sgn (\pi)\tilde\omega(\theta_{\pi (1)},\theta_{\pi (2)})\dots \tilde\omega(\theta_{\pi (2k-1)},\theta_{\pi (2k)})\omega^k
\label{wedging}
\end{equation}
Thus wedging together the $2k$ free $dz_a$ indices can again be evaluated by contractions with $\tilde \omega$, so this characteristic number is indeed a Rozansky-Witten invariant. It is obtained (up to some overall factor) by evaluating a sum of graphs, each one of which is obtained from a permutation $\pi$ of the trivalent vertices of the wheel: we join together the spokes from $\pi(1)$ and $\pi(2)$, the spokes $\pi(3)$ and $\pi(4)$ etc. We must keep track of signs however, and so we have to consider the orientations more closely. 
\vskip .25cm
 A trivalent graph in the plane acquires a canonical orientation by taking the anticlockwise cyclic ordering of the edges at each vertex. To connect this with the orientation needed for the Rozansky-Witten invariant we follow essentially the approach of Kapranov \cite{Kap} and introduce the notion of a {\it flag}: an edge together with a choice of a vertex lying on it. Note that a flag is the same thing as an oriented edge: we choose the distinguished vertex to be the initial point of the arrow. For any finite set $S$, let $\det S$ denote the highest exterior power of the vector space $\R^{S}$ of maps from $S$ to $\R$. Then the planar orientation is the same as an orientation on the one-dimensional vector space
$$\bigotimes_{v\in V(\Gamma)} \det F(v)$$
where $V(\Gamma)$ is the set of vertices of $\Gamma$ and $F(v)$ the three-element set of flags whose distinguished vertex is $v$. Now let $F(\Gamma)$ be the space of all flags in $\Gamma$, then 
$\det F(\Gamma)$
is the top exterior power of the vector space
$$\bigoplus_{f\in F(\Gamma)}{\R}^{\{f\}}=\bigoplus_{v\in V(\Gamma)}{\R}^{F(v)}$$
The right hand side is a direct sum of three-dimensional spaces parametrized by the vertices of $\Gamma$. Since three-forms anticommute, we have
\begin{equation}
\det F(\Gamma)\cong \det V(\Gamma)\otimes \left(\bigotimes_{v\in V(\Gamma)} \det F(v)\right)
\label{orient}
\end{equation}
On the other hand, each edge occurs twice in the set of flags, since it has two ends, so we have
$$\bigoplus_{f\in F(\Gamma)}{\R}^{\{f\}}=\bigoplus_{e\in E(\Gamma)}{\R}^{F(e)}$$
where $F(e)$ is the two-element set of flags containing the edge $e$. Now two-forms commute, so 
\begin{equation}
\det F(\Gamma)\cong \bigotimes_{e\in E(\Gamma)} \det F(e)
\label{orient1}
\end{equation}
An orientation of the two dimensional space ${\R}^{F(e)}$ is just given by choosing which basis element to put first, and this is a choice of flag in $F(e)$ or, as we saw, a choice of orientation on the edge $e$. We defined the Rozansky-Witten orientation as an equivalence class of orderings 
  of vertices   together with orientations on the  edges, so in our current formalism this is an orientation on 
$$\det V(\Gamma) \otimes \left(\bigotimes_{e\in E(\Gamma)} \det F(e)\right)$$
But from (\ref{orient}) and (\ref{orient1}) this is equivalent to an orientation given by cyclic orderings at the vertices.
\vskip .25cm
 Let us apply this to the situation of the wheel. Numbering the vertices $\{1,2\dots,2k\}$, the contractions around the hub of the wheel correspond to giving the edges on the hub the clockwise orientation. Let the flag given by the edge to the right of vertex $1$ be denoted $f_1$ etc. and call the flag corresponding to the opposite orientation $\bar f_1$. Call the flag emanating along a spoke from $1$, $\varphi_1$ etc. Give $\det V(\Gamma)$ the orientation $v_1\wedge \dots \wedge v_{2k}$, then the planar orientation is determined by the following element in $\det F(\Gamma)$: 

$$(\bar f_{2k}\wedge f_1 \wedge \varphi_1)\wedge (\bar f_{1}\wedge f_2 \wedge \varphi_2)\wedge \dots (\bar f_{2k-1}\wedge f_{2k} \wedge \varphi_{2k})=$$ 
\begin{equation}-(f_1 \wedge \bar f_{1})\wedge (f_2 \wedge \bar f_{2})\wedge \dots\wedge (f_{2k}\wedge \bar f_{2k})\wedge (\varphi_1\wedge \varphi_2\wedge \dots\wedge \varphi_{2k})
\label{sign2}
\end{equation}
Now if, for a permutation $\pi \in S_{2k}$, we pair together in an ordered way the spokes $(\pi(1),\pi(2)),(\pi(3),\pi(4)),\dots$ the orientation on the space spanned by the corresponding flags is 
$$(\varphi_{\pi(1)}\wedge \varphi_{\pi(2)})\wedge \dots\wedge (\varphi_{\pi(2k-1)}\wedge\varphi_{\pi(2k)})=(\sgn \pi) \varphi_1\wedge \varphi_2\wedge \dots\wedge \varphi_{2k}$$
Thus, using (\ref{sign2}), in the planar orientation, each summand in (\ref{wedging}) gives a contribution $-1$ for the unordered pairing of spokes defined by $\pi$.
\vskip .25cm
Finally, using (\ref{sign1}), (\ref{sign2}) and (\ref{wedging}) we have a graphical formula for the characteristic number $s_{2k}$:
$$s_{2k}[M]=\int_M \Gamma(\Phi)\omega^k$$
where $\Gamma$ is the linear combination of graphs (with the planar orientation) 
$$\Gamma=(-1)^k (-1) \frac{1}{(2\pi i)^{2k}} \frac{1}{2^{2k}(k!)^2} 2^k k!\sum \mathrm{pairings \,\, of\,\, spokes\,\, of\,\, a\,\,} 2k {\mathrm{-wheel}}$$
The coefficient simplifies to 
\begin{equation}
-\frac{1}{(8\pi^2)^k k!}
\label{coeff}
\end{equation}
\vskip .25cm
The process of connecting in pairs the spokes of a wheel was essentially a graphical interpretation of the exterior product, and we can use the same formalism to deal with characteristic classes $s_{2k_1}s_{2k_2}\dots s_{2k_m}$
where $k_1+k_2+\dots + k_m=k$. These correspond up to a scalar multiple to the differential form
$$(\tr K^{2k_1})\wedge (\tr K^{2k_2})\wedge \dots \wedge(\tr K^{2k_m})$$
and this is graphically obtained by taking the disjoint union of wheels $w_{2k_1},\dots, w_{2k_m}$, summing over all ways of pairing the $2k$ spokes and multiplying by a  factor.  As in (\ref{sign2}) we obtain a factor $-1$, but now one for each wheel, giving an overall factor
\begin{equation}
(-1)^m\frac{1}{(8\pi^2)^k k!}
\label{coeff1}
\end{equation}
 Now any characteristic number is a linear combination of such products since any symmetric polynomial is a polynomial in the sums of powers, so we see that characteristic numbers as Rozansky-Witten invariants arise from the graphs generated by wheels. Fortunately these have been studied in some detail by three-manifold theorists.

\section{Wheeling} 
We refer to \cite{BGRT} for details of the following facts in graph homology. For us a {\it unitrivalent} graph is a possibly disconnected graph whose vertices are either univalent or trivalent. A disjoint union of wheels and purely trivalent graphs is an example.
Let ${\cal B}'$ be the space of linear combinations of unitrivalent graphs  modulo anti-symmetry and the IHX relation on the trivalent vertices. There are two products on ${\cal B}'$. The first, denoted $\cup$, is induced by disjoint union  of graphs and the second, denoted $\times$,  arises from a vector space isomorphism of ${\cal B}'$ with the algebra of chord diagrams. We shall only be using part of this structure, and indeed not the full force of the Wheeling Theorem, and for us it is sufficient to remark that the product $\times$ adds the number of univalent vertices and so the purely trivalent graphs form a subalgebra, and for this subalgebra  the multiplication $\times$ is again induced by disjoint union. Both algebra structures are graded by half the number of vertices. 

Given a unitrivalent graph $C$, we obtain an operator $\hat{C}:{\cal
B}^{\prime}\rightarrow{\cal B}^{\prime}$ defined in the following
way. If $C$ has no more univalent vertices than $C^{\prime}$, then
$\hat{C}(\Gamma)$ is defined by summing over all the ways of joining
them to  the univalent vertices of $\Gamma$;
otherwise we define it to be zero. We  then extend this linearly to any element
$C\in{\cal B}^{\prime}$. 
\vskip .25cm
Since this is close to what we have been doing, let us look at a particular example, taking  
$$C=w_{2k},\qquad \Gamma ={\ell}^k$$
where $\ell$ is the line: the unique graph with one edge and two vertices. In this case we calculate $\hat C(\Gamma)$ by summing over all ways of joining the $2k$ ends of the $k$ disjoint lines to the $2k$ univalent vertices of the wheel $w_{2k}$. This is $2^k k!$ times the sum of all ways of joining up the spokes of the wheel in pairs. Denote this latter summing process by $C \mapsto S(C)$.

The Wheeling Theorem concerns the special role of the following element in ${\cal B}'$:
\begin{equation}
  \Omega=\exp_{\cup}  \sum_{n=1}^\infty b_{2n}w_{2n}.
  \label{O}
\end{equation}
where
\begin{equation} 
  \sum_{n=0}^\infty b_{2n}x^{2n} = \frac{1}{2}\log\frac{\sinh x/2}{x/2}.
  \label{sinh}
\end{equation}
Thus
\begin{equation}
\Omega=1+\frac{1}{48}w_2+\frac{1}{2!48^2}(w_2^2-\frac{4}{5}w_4)+\dots
\label{expand}
\end{equation}
Here $\exp_{\cup}$ means that we take products using disjoint union of graphs. The theorem, conjectured in \cite{BGRT}, and independently by Deligne \cite{D},  states:
\begin{thm}
The operator associated to $\Omega$ intertwines the two product
structures on ${\cal B}^{\prime}$; ie.\ $\hat{\Omega}:{\cal
B}^{\prime}_{\cup}\rightarrow{\cal B}^{\prime}_{\times}$ is an
isomorphism of algebras.
\end{thm}
This was recently proved by  Bar-Natan, Le, and
Thurston~\cite{blt}. It is also a corollary of Kontsevich's results on deformation quantization \cite{Kon}.
\vskip .25cm
A particular case of  the Wheeling Theorem is
\begin{equation}
\hat \Omega (\ell^k_{\cup})=(\hat \Omega(\ell))^k_{\times}
\label{wh}
\end{equation}

We remarked above that applying $\hat C$ to a product of lines involves summing over the pairings of free edges and multiplying by $2^k k!$, so using this and (\ref{expand}) the right hand side of (\ref{wh}) is
$$(\ell + \frac{1}{24}\Theta)^k_{\times}$$
whose purely trivalent part  is
\begin{equation}
(\frac{1}{24}\Theta)^k_{\cup}
\label{RHS}
\end{equation}
Now consider the left hand side of (\ref{wh}). Using the summing process $S$ we obtain  Rozansky-Witten invariants
$$\int_M S(w_{2k_1}\cup \dots \cup w_{2k_m})(\Phi)\omega^k=(-1)^m {(8\pi^2)^k k!} (s_{2k_1}\wedge \dots \wedge s_{2k_m})[M]$$
This means that
$$\int_M S(\exp_{\cup}  \sum_{n=1}^\infty b_{2n}w_{2n})(\Phi)\omega^k= (8\pi^2)^k k!\exp_{\wedge}  ( -\sum_{n=1}^\infty b_{2n}s_{2n}) [M]$$
But each characteristic class $s_{2n}$ is $\sum x_i^{2n}$ for the Chern roots $x_i$, so
\begin{equation}
\exp ( -\sum_{n=1}^\infty b_{2n}s_{2n})[M]=\prod_i \exp (-\sum_{n=1}^\infty b_{2n}x_i^{2n})[M]=\prod_i \left(\frac {x_i/2}{\sinh x_i/2}\right)^{1/2}[M]
\label{ex}
\end{equation}
from (\ref{sinh}). Now since $\sum x_i=0$ 
$$\Td=\prod_i\frac{x_i}{1-e^{-x_i}}=\prod_i\exp(x_i/2)\frac{x_i/2}{\sinh x_i/2}=\prod_i\frac{x_i/2}{\sinh x_i/2}$$
so the right hand side of (\ref{ex}) is  ${\Td}^{1/2}={\hat A}^{1/2}$. Equating the trivalent parts of the left and right hand sides of (\ref{wh}), and evaluating on $[M]$ we obtain
$$
\frac{1}{(24^k)}(8\pi^{2})^k k!b_{\Theta^k}(M)=2^k k!{(8\pi^2)^k k!}\sqrt{\hat A}\,[M]$$
Finally, using the formula  (\ref{bb}) for $b_{\Theta^k}$, we find
\begin{thm} Let $M$ be a compact irreducible hyperk\"ahler manifold, then 
$$\frac{1}{(192\pi^2 k)^k}  \frac{\Vert R \Vert ^{2k}}{(\vol M)^{k-1}}=\sqrt{\hat A}\,[M]$$
\end{thm}

\section {An inequality}

An immediate consequence of the formula in Theorem 3 is that if $M$ is a compact irreducible hyperk\"ahler manifold, then
$$\sqrt{\hat A}\,[M]>0.$$
Equality would mean that $\Vert R \Vert=0$, implying that  $M$ is flat and thus contradicting irreducibility. We can add to this inequality the following consequence of Proposition 1:
$${\hat A}[M^{4k}]=\Td[M]=\sum(-1)^p h^{0,p}=k+1$$
\vskip .25cm
For $k=1$, the K3 surface, we have 
$$\sqrt{\hat A}\,[M]=\frac{1}{2}{\hat A}[M]=1$$
\vskip .25cm
For $k=2$,
$$\sqrt{\hat A}\,[M]=\frac{1}{2}\hat A_2 [M]-\frac{1}{8} \hat A_1^2[M]=\frac{3}{2}-\frac{1}{8}\hat A_1^2[M]$$
so our inequality reads
$$\hat A_1^2[M]<12$$
Using
$$\hat A_1=\frac{1}{12} c_2,\qquad \hat A_2=\frac{1}{720}(-c_4+3c_2^2)$$
this gives an inequality for the Euler characteristic $\chi(M)=c_4[M]$ which is
$$\chi(M)<3024$$
Although these numbers have high divisibility in general (see \cite{Sal1}) this is still rather a crude estimate in this dimension. Using the Riemann-Roch theorem and the Bogomolov/Verbitsky results on the cohomology ring generated by $H^2(M)$, Beauville \cite{Bea1} has given a sharp estimate for $k=2$:
$$\chi(M)\le 324.$$

\end{document}